\documentclass[12pt]{amsart}
\usepackage{amssymb,amsmath,amsthm,MnSymbol}
\usepackage[linktocpage=true,backref=page]{hyperref}
\hypersetup{
	colorlinks=true, 	%set true if you want colored links
	linktoc=all,     		%set to all if you want both sections and subsections linked
	linkcolor=magenta,  	%choose some color if you want links to stand out
	citecolor=cyan,
}
\calclayout
\newtheorem{theorem}{Theorem}[section]
\numberwithin{theorem}{section} \numberwithin{equation}{section}
\title[Holomorphic Conformal Structures and K3 Surfaces]{On holomorphic conformal structures associated with lattice polarized K3 surfaces}
\author{Andreas Malmendier}
\address{Dept.\!~of Mathematics \& Statistics, Utah State University, Logan, UT 84322}
\email{andreas.malmendier@usu.edu}
\thanks{A.M. acknowledges support from the Simons Foundation through grant no.~202367.}
\author{Michael T.  Schultz}
\address{Dept.\!~of Mathematics, Virginia Tech, Blacksburg, VA 24060}
\email{michaelschultz@vt.edu}
\begin{document}
\begin{abstract}
We discuss the connection between Picard-Fuchs equations for certain families of lattice polarized K3 surfaces and the construction of integrable holomorphic conformal structures on their period domains. We then compute an explicit example of a locally conformally flat holomorphic metric associated with generic Jacobian Kummer surfaces, which allows for a novel description of the local variation of complex structure. 
 \end{abstract}
\maketitle
\section{introduction}
\label{s-intro}
Among complex manifolds, the existence of holomorphic conformal structures are quite restrictive \cite{MR0685426}, as opposed to the well known conformal structures in the real Riemannian or semi-Riemannian case \cite{MR0355886}. In this article, we will explore holomorphic conformal structures that are integrable, or equivalently, locally conformally flat. Quadric hypersurfaces $Q^n \subset \mathbb{P}^{n+1}$ naturally possess such an integrable holomorphic conformal structure and serve as the model spaces of (holomorphic) conformal geometry; see Kuiper \cite{MR0031310} and Kobayashi \cite[IV. Theorem 6.2]{MR0355886}.  As we will explain, the study of integrable holomorphic conformal structures, well known in the field of differential geometry,  finds another powerful application in the field of algebraic geometry, namely within the study of algebraic families of certain algebraic varieties. 
\par Such a situation arises when a proper smooth  holomorphic map $f \, \colon \, \mathfrak{X} \to M$ is given between complex algebraic varieties, where $\mathfrak{X}, M$, and the general fiber $X_t=f^{-1}(t)$ for $t \in M$ are smooth and connected, and a distinguished class of projective embeddings $\mathfrak{X} \subset \mathbb{P}^N$ is given; see \cite{MR258824}. Fix a generic point $t=0 \in M$. Then we can think of this situation as a family  $\{ X_t \}_{t \in M}$ of generically smooth, compact, complex projective algebraic varieties with the algebraic parameter space $M$. However, in general $M$ is not  complete, and $X_t$ acquires singularities as $t$ tends to $\infty$ or, more generally, towards a more complicated discriminant locus contained in $M$ with $t \not =0$.  The fundamental problem is to describe the variation of the complex structure of the compact algebraic manifold $X_0$.
\par In this article, we will investigate algebraic families of varieties arising as lattice polarized K3 surfaces.  A family of lattice polarized K3 surfaces can be constructed for every even, primitive sublattice $L$ of rank $\rho =20-n > 0$ of the K3 lattice $\Lambda_{\mathrm{K3}} \cong 3 H \oplus 2 E_8(-1)$ of signature $(3,19)$. These $L$-polarized K3 surfaces are classified, up to isomorphism, by a coarse moduli space $\mathfrak{M}_L$, which is a quasi-projective variety of dimension $n$.  Given an algebraic family of $L$-polarized K3 surfaces $f_L \,\colon\, \mathfrak{X} \to M$ over an $n$-dimensional open subspace $M$, disjoint from the discriminant locus, one can compute a holonomic system $E$ of rank $n+2$ in $n$ variables, whose  solutions are the period integrals of $X_t$. This linear system of partial differential equations is called the \emph{Picard-Fuchs system} for the family $f_L$.   This system is a fundamental object of study, encoding crucial information about the geometry of the K3 surface as the complex structure varies with $t \in M$. In fact, a connection on a vector bundle over the base space $M$ can now be constructed where the fibers are the direct sum of the Hodge groups $H^{p, q}(X_t)$. As such a family is locally topologically trivial, cohomology classes can be moved from one fiber in the family to nearby fibers, providing flat local sections of the bundle. In turn, from these flat sections, the existence of such a connection can be inferred, called the Gauss-Manin connection \cite{MR229641}.  
\par The goal of this article is to explain the great utility of the framework of integrable holomorphic conformal structures in the study of Picard-Fuchs equations for families of lattice polarized K3 surfaces. While this connection has been known for quite sometime, it does not appear to have been used much since the pioneering work of Sasaki \& Yoshida  \cite{MR0996019}.  We will demonstrate how the Picard-Fuchs system for a family $f_L \, \colon \, \mathfrak{X} \to M$ of $L$-polarized K3 surfaces imparts $M$ with an integrable holomorphic conformal structure. In fact, the period map for $f_L$ determines a (multivalued) inverse of the canonical projection $\mathcal{D}_L \to \mathfrak{M}_L$ where $\mathcal{D}_L$ denotes the period domain. The latter can be identified with a relatively open subset of a quadric hypersurface in $\mathbb{P}^{n+1}$, and thus admits a natural integrable conformal structure. While the holomorphic conformal structure is sufficiently interesting in its own right, we will see that in this context it efficiently encodes much more information -- essentially, the entire local variation of Hodge (or equivalently in this case, complex) structure. 
\par This article is structured as follows: after reviewing some aspects of lattice polarized K3 and Kummer surfaces and the geometry of their period domains, we provide a new example of an integrable holomorphic conformal structure associated with Jacobian Kummer surfaces, allowing for a novel description of their local variation of complex structure. 
\section{Lattice Polarized K3 \& Kummer Surfaces}
\label{s-k3_kummer}
Let us start by reviewing a few classical lattice theory facts. We shall use the following standard notations: $L_1 \oplus L_2$ is orthogonal direct sum of the two lattices $L_1$ and $L_2$, $nL = L^{\oplus n}$ is the $n$-fold orthogonal direct sum of $L$, $L(\lambda)$ is obtained from the lattice $L$ by multiplication of its form by $\lambda \in \mathbb{Z}$, $\langle R \rangle$ is a lattice with the matrix $R$ in some basis; $A_n$, $D_m$, and $E_k$ are the positive definite root lattices for the corresponding root systems;  $H$ is the unique even unimodular hyperbolic rank-two lattice.  Given a lattice $L$, one has the discriminant group $A(L) = L^\vee/L$ and its associated discriminant form, denoted $q_L$.  For a lattice $L$ and a field $\mathbb{K}$ of characteristic zero, let $L_{\mathbb{K}} = L \otimes_{\mathbb{Z}} \mathbb{K}$. We will continue to denote by $( \, \cdot \, , \, \cdot \, )$, by a slight misuse of notation,  the $\mathbb{K}$-extension of the bilinear form on $L$ to $L_{\mathbb{K}}$.
	
\par A K3 surface $X$ is a smooth compact complex surface with trivial canonical bundle $\omega_X := \bigwedge^2  \Omega^1_X \cong \mathcal{O}_X$ and $h^1(X, \mathcal{O}_X) = 0$. Here, $\Omega^1_X$ is the holomorphic cotangent bundle of $X$. It is known that every K3 surface is K\"ahler.  The second cohomology group $H^2(X, \mathbb{Z})$, equipped with the pairing induced by the cup product $(\,\cdot\, ,\,\cdot\,)$, forms a lattice isometric to the even rank-22 K3 lattice $\Lambda_{\mathrm{K3}} \cong 3 H \oplus 2 E_8(-1)$ of signature $(3,19)$. Up to isometry, $\Lambda_{\mathrm{K3}}$ is the only such even lattice with the given signature. The isometry is called the marking of the K3 surface $X$. Thus, we have
\begin{equation*}
\Lambda_{\mathrm{K3}} \otimes \mathbb{C} \ \cong \ H^2(X,\mathbb{C}) = H^{2,0}(X) \oplus H^{1,1}(X) \oplus H^{0,2}(X)
\end{equation*}
by the Hodge decomposition.  The N\'eron-Severi lattice is given by
\begin{equation}
 \mathrm{NS}(X) := H^{1,1}(X) \cap H^2(X, \mathbb{Z}).
\end{equation}
Via the first Chern class map it is  isomorphic to the Picard lattice $\mathrm{Pic}(X)$. This lattice is even and has signature $(1, \rho_X-1)$, where $1 \leq \rho_X \leq 20$ is the Picard number of $X$.  The transcendental lattice $T_X$ is defined to be the smallest sublattice of $H^2(X, \mathbb{Z})$  whose complexiﬁcation contains a generator of $H^{2,0}(X)$. In the non-degenerate case, it follows that $T_X$ is the orthogonal complement of $\mathrm{NS}(X)$ in $H^2(X, \mathbb{Z})$ and has signature $(2, 20 - \rho_X)$. Conversely, for every primitive sublattice $L$ (or $N$) of $\Lambda_{\mathrm{K3}}$ of signature $(1, \rho_X-1)$ (or $(2, 20 - \rho_X)$), there is an algebraic K3 surface $X$ together with a marking such that $L \cong \mathrm{NS}(X)$ (or $N \cong T_X$). Here, the term algebraic surface refers to a complex algebraic variety of dimension two.  
\par A (pseudo-)polarized K3 surface of degree $2k$ with $k>0$ is a pair $(X, h)$ consisting of a K3 surface $X$ and a primitive (pseudo-)ample class $h \in  \mathrm{NS}(X)$ with $( h, h ) =2k$. Moreover, two (pseudo-)polarized K3 surfaces $(X, h)$ and $(X', h')$ of degree $2k$ are equivalent if there exists an analytic isomorphism $\alpha \, \colon \, X \to X'$ such that $\alpha^* h' =h$. If $\mathrm{NS}(X)$ contains a pseudo-ample class, then it follows that $X$ is projective, hence every pseudo-polarized K3 surface is projective \cite{MR3586372}.
\par There are important classes of examples of K3 surfaces that arise as the result of a pseudo-polarization. The first is the class of double sextic surfaces; they are pseudo-polarized K3 surfaces $(X, h)$ of degree 2 with an effective divisor $D$ on $X$ such that $[D]=h \in \mathrm{NS}(X)$: in the generic case, the linear system $|D|$ then defines a map from $X$ to $\mathbb{P}^2$ that is generically $2 \colon 1$. Accordingly, $X$ is birational to the double cover of $\mathbb{P}^2$ branched on a sextic curve.   A second class of examples are quartic hypersurfaces. They are  pseudo-polarized K3 surfaces $(X, h)$ of degree 4 with an effective divisor $D$ on $X$ such that $[D]=h$: in the generic case, $|D|$ defines a birational morphism onto a quartic hypersurface in $\mathbb{P}^3$. Quartic hypersurfaces, for example the Fermat quartic or the Hudson quartic, constitute some of the famous examples of complex projective K3 surfaces. For degree 4, there is also a second situation that can arise, called a double quadric surface. Here, $|D|$ defines a generically $2 \colon 1$ map from $X$ onto a quadric hypersurface in $\mathbb{P}^3$ that is isomorphic to $\mathbb{P}^1 \times \mathbb{P}^1$.  Accordingly, $X$ is birational to the double cover of $\mathbb{P}^1 \times \mathbb{P}^1$ branched on a curve of bi-degree $(4, 4)$.

\par Let $X$ be a smooth complex projective K3 surface.  A lattice polarization on $X$ is defined as a primitive lattice embedding $\jmath \, \colon \, L \hookrightarrow \mathrm{NS}(X)$, where $\jmath(L)$ contains a pseudo-ample class. Here, $L$ is an even indefinite lattice of signature $(1,\rho_L-1)$, with $1 \leq \rho_L \leq 20$. Two K3 surfaces $(X,\jmath)$ and $(X',\jmath')$ are said to be isomorphic under $L$-polarization if there exists an analytic isomorphism $\alpha \colon X \rightarrow X'$ and a lattice isometry $\beta \in O(L)$ such that $\alpha^* \circ \jmath' = \jmath \circ \beta$, where $\alpha^*$ is the induced morphism at the cohomology level.\footnote{Our definition of isomorphic lattice polarizations coincides with the one used by Vinberg, and is slightly more general than the one used in \cite[Sec.~1]{MR1420220}.} $L$-polarized K3 surfaces are classified up to isomorphism by a coarse moduli space $\mathfrak{M}_L$, which is a quasi-projective variety of dimension $20-\rho_L$, see \S \ref{ss-periods}. A general $L$-polarized K3 surface $(X, \jmath)$ satisfies $\jmath(L)=\mathrm{NS}(X)$.
\subsection{Kummer surfaces}
\label{ss:Kummers}
We also consider special K3 surfaces known as Kummer surfaces; they are associated with abelian surfaces. First, we recall that a complex 2-dimensional torus $T^2$  is a surface isomorphic to the quotient of $\mathbb{C}^2$ by a real non-degenerate lattice $G$ of rank 4. The group $H^2(T^2, \mathbb{Z})$,  equipped with the pairing induced by the cup product, is an even unimodular lattice of signature $(3, 3)$. Up to isometry, there exists only one lattice with these properties, namely the lattice $3H$. An abelian surface $A$ is an algebraic torus.  Moreover, for every primitive sublattice $S$ (or $T$) of $3H$ of signature $(1, \rho-1)$ (or $(2, 4 -\rho)$) with $1 \le \rho \le 4$, there is an abelian surface $A$ with a marking such that $S \cong \mathrm{NS}(A)$ (or $T \cong T_A$).  Every abelian surface admits an involution $\imath$: the involution $\imath$ maps every element $a \in A$ to its inverse $-a \in A$ with respect to the group law. The fixed points of this map are the sixteen points of order two on $A$.  The Kummer surface $Y$ associated with $A$ is the K3 surface obtained as the minimal resolution of the quotient $A/\langle \imath \rangle$. Let $\widetilde{A}$ be the blow up of $A$ in the sixteen fixed points of the involution $\imath$ and let $\pi \colon \widetilde{A} \to Y$ be the $2 \, \colon \, 1$ cover between $\widetilde{A}$ and $Y$.  One can easily check that $\pi_* H^2(\widetilde{A}, \mathbb{Z})^\imath  \cong 3 H(2)$.

\par By construction, the Kummer surface $Y$ admits an even set of sixteen rational curves, coming from the blow up of the order-2 points on the abelian surface, with classes $K_1, \dots, K_{16}$. The term even refers to the fact that the classes add up to an even element in $\mathrm{NS}(Y)$, that is, $K_1 + \dots + K_{16} \in 2 \mathrm{NS}(Y)$.  An even set made up of disjoint rational curves on any K3 surface must contain either eight or sixteen curves \cite{MR0429917}.  The classes  $K_1, \dots, K_{16}$ on $Y$ generate a minimal primitive sublattice $K$ of $\mathrm{NS}(Y)$, and hence a sublattice of $\Lambda_{\mathrm{K3}}$. The lattice $K$ is the orthogonal complement of $\pi_* H^2(\widetilde{A}, \mathbb{Z})^\imath$ in $\Lambda_{\mathrm{K3}}$, whence a negative definite even lattice of rank 16 and determinant of the discriminant form $2^6$. It was described explicitly in \cite{MR0429917}.

\par A polarization of an abelian surface $A \cong \mathbb{C}^2/G$, over the complex numbers, is a  positive definite hermitian form $H$ on $\mathbb{C}^2$,  satisfying $\operatorname{Im} H(\Lambda,\Lambda) \subset \mathbb{Z}$.  Every such hermitian form determines a line bundle $\mathcal{L} \to A$ in the N\'eron-Severi lattice $\mathrm{NS}(A)$.  A general fact from linear algebra asserts that one can always  choose a basis of $G$ such that $E = \operatorname{Im} H$ is given  by the matrix $\bigl(\begin{smallmatrix} 0&M\\ -M&0 \end{smallmatrix} \bigr)$ with $M=\bigl(\begin{smallmatrix}d_1&0\\ 0&d_2 \end{smallmatrix} \bigr)$ where $d_1, d_2 \in \mathbb{N}$, $d_1, d_2 \ge 0 $, and $d_1$ divides $d_2$. The pair $(d_1, d_2)$ gives the type of the polarization; see \cite{MR2062673}.  A polarization of type $(1,1)$ is also called a principal polarization. 

\par For a principally polarized abelian variety $A$ the line bundle $\mathcal{L}$ defining its  principal polarization is ample and satisfies $h^0(\mathcal{L}) = 1$. There exists an effective divisor $\Theta$ such that $\mathcal{L}= \mathcal{O}_A(\Theta)$, uniquely defined only up to translations. The divisor $\Theta \in \mathrm{NS}(A)$ is called a theta divisor associated with the polarization.  Moreover, it is known that a principally polarized abelian surface is either the Jacobian of a smooth curve $C$ of genus 2 with canonical theta divisor \cite[\S 11.1]{MR2062673}, or the product of two complex elliptic curves, with the product polarization \cite{MR2514037}.  
\par For the principally polarized Abelian surfaces $A = \mathrm{Jac} \,C$ where $C$ is a smooth curve of genus 2 with $T_A = 2H + \langle -2 \rangle$, the general Kummer surface associated with $A$ satisfies 
\begin{equation}
\label{eqn:lattice}
 \mathrm{NS}(Y) \ = \ H \oplus D_8(-1) \oplus D_4(-1) \oplus A_3(-1)  \cong \ H \oplus D_7(-1) \oplus 2 D_4(-1) \,.
\end{equation} 
This follows from the results in \cite{MR3263663}.  We also refer to the K3 surface $Y$ as a \emph{generic} Jacobian Kummer surface. A similar statement holds for Kummer surfaces associated with  the product of two non-isogenous elliptic curves such that $T_A = 2H$, a specialization of the Kummer surface where the lattice polarization extends. See ~\cite{MR2409557}.
\subsection{Double sextic surfaces}
\label{ss:double_sextics}
The Jacobian Kummer surfaces $Y$ can be described explicitly as double sextic surfaces as follows:  let $C$ be given in affine coordinates $(\xi,\eta)$ by the Rosenhain normal form 
\begin{equation}
\label{eqn:Rosenhain}
 C\colon \quad \eta^2 = \xi \,\big(\xi-1) \, \big(\xi- \lambda_1\big) \,  \big(\xi- \lambda_2 \big) \,  \big(\xi- \lambda_3\big) \,.
\end{equation} 
We denote the hyperelliptic involution on $C$ by $\imath_C$. An ordered tuple $(\lambda_1, \lambda_2, \lambda_3)$ -- where the $\lambda_i$ are pairwise distinct and different from $(\lambda_4,\lambda_5,\lambda_6)=(0, 1, \infty)$ -- determines a point in the moduli space of genus-2 curves with level-two structure.\footnote{A full level-2 structure is basis of rational 2-torsion points on the principally polarized abelian surface $\mathrm{Jac}\, C$.}   
\par For the symmetric product $C^{(2)}$, the quotient $C^{(2)}/\langle \imath_C \times  \imath_C \rangle$ is realized as a variety in terms of variables $t=\xi^{(1)}\xi^{(2)}$, $x=\xi^{(1)}+\xi^{(2)}$, and $y=\eta^{(1)}\eta^{(2)}$, satisfying the affine equation
\begin{equation}
\label{kummer_middle}
  y^2 = t \big(  t  - x +  1 \big)  \prod_{i=1}^3 \big( \lambda_i^2 \, t  -  \lambda_i \, x +  1 \big) \,.
\end{equation}
The affine variety completes to a hypersurface in $\mathbb{P}(1,1,1,3)$ which is birational to the Kummer surface $\operatorname{Kum}(\operatorname{Jac} C)$, i.e., we have 
\begin{equation}
 Y \ \cong \ \operatorname{Kum}(\operatorname{Jac} C) \,.
\end{equation}
We introduce the following six lines $\ell_1, \dots, \ell_6$ in $\mathbb{P}^2 = \mathbb{P}(z_1, z_2, z_3)$:
\begin{equation}
\label{eqn:lines}
  \ell_1= z_1 \,, \ \ell_2 = z_2\,, \ \ell_3 = z_1 + z_2 - z_3 \,, \  \ell_{3+i} = z_1 + \lambda_i^2 z_2 - \lambda_i z_3 \ (\text{for $i=1, 2, 3$}) \,,
\end{equation}  
which are tangent to a common conic $z_3^2 -4 z_1 z_2=0$. Then, Equation~\eqref{kummer_middle} is the double cover of $\mathbb{P}^2$ branched on $\ell_1, \dots, \ell_6$ in the affine chart $z_1=1, z_2=t, z_3=x$. 
\par Generalizing further, we consider the minimal resolution of a double cover of the projective plane $\mathbb{P}^2$ branched along the union of six lines in general position. We denote the six-line configuration by $\boldsymbol \ell =\{\ell_1,\dots,\ell_6\}$. Using the weighted homogeneous coordinates $[z_1 : z_2 : z_3 : y] \in \mathbb{P}(1,1,1,3)$, the family of double sextics is given by 
\begin{equation}
\label{eqn-six_lines}
 	y^2=\prod_{i=1}^6\big(a_{i1} z_1+a_{i2} z_2+a_{i3} z_3\big),
\end{equation}
where the lines $\ell_i =\left\{ [z_1 : z_2 : z_3] \; | \; a_{i1} z_1+a_{i2} z_2+a_{i3} z_3 = 0 \right\} \subset \mathbb{P}^2$ for parameters $a_{i j} \in \mathbb{C}$, $i=1,\dots,6$, $j=1,2,3$ are assumed to be in general position.  Let $ (a_{i j}) \in \mathrm{Mat}(3,6;\mathbb{C})$ be the matrix whose entries are the coefficients encoding the six-line configuration $\boldsymbol \ell$. The general member of the family in Equation~(\ref{eqn-six_lines}) is  a smooth K3 surface of Picard number 16. It follows from \cite{MR4015343} that Equation~(\ref{eqn-six_lines}) defines a family of $L$-polarized K3 surfaces, where  $L$ has rank 16 with the following isomorphic presentations:
\begin{equation}
\label{eqn:L}
\begin{split}
 L  \ \cong \ & H \oplus E_8(-1)  \oplus 6 A_1(-1)  \ \cong  \  H \oplus D_8(-1) \oplus D_4(-1) \oplus 2 A_1(-1).
\end{split}
\end{equation}
\par Let $\mathcal{D}$ be the configuration space of six lines $\boldsymbol \ell$ whose minimal resolution is a K3 surface. Note that isomorphic K3 surfaces are obtained if we act on an element $(a_{i j}) $ by matrices induced by automorphisms of $\mathbb{P}^2$ on the left and overall scale changes of each line $\ell_i \in \boldsymbol \ell$ on the right. Thus, we are led to consider the four-dimensional quotient space
\begin{equation}
\label{eqn:M6}
\mathcal{M} = \mathrm{SL}(3,\mathbb{C}) \backslash \mathcal{D} / (\mathbb{C}^*)^6.
\end{equation}
\par The family of double-sextics in Equation~(\ref{eqn-six_lines}) has been studied in the literature, famously by Matsumoto \cite{MR1103969},  and Matsumoto et al. \cite{MR973860, MR1073363, MR1136204}. 
In \cite{MR973860},  the authors used the action of $\mathrm{SL}(3,\mathbb{C})$ and $(\mathbb{C}^*)^6$ to bring the matrix $(a_{i j})$  into the form
\begin{equation}
\label{eqn-matsumoto_sextic}
\begin{pmatrix}
1 & 0 & 0 & 1 & 1 & 1 \\
0 & 1 & 0 & 1 & x_1 & x_2 \\
0 & 0 & 1 & 1 & x_3 & x_4 
\end{pmatrix},
\end{equation}
so that the associated K3 surfaces are the minimal resolution of the double sextic surfaces
\begin{equation}
\label{eqn:MATeqn}
 \mathcal{X}\colon \quad y^2 = z_1 z_2 z_3 \big(z_1 + z_2 + z_3\big) \big(z_1 + x_1 z_2 + x_3 z_3\big)  \big(z_1 + x_2 z_2 + x_4 z_3\big) \,.
\end{equation}
\par A six-line configuration $\boldsymbol \ell$ is tangent to a conic if and only if a certain projective invariant $R$ vanishes; see \cite{MR1007155}. For a six-line configuration given by Equation~(\ref{eqn-matsumoto_sextic}), one finds
\begin{equation}
\label{eqn:R}
 R =
-x_{1} x_{2} x_{3}+x_{1} x_{2} x_{4}+x_{1} x_{3} x_{4}-x_{2} x_{3} x_{4}-x_{1} x_{4}+x_{2} x_{3} .
\end{equation}
The subvariety determined by $R=0$ is also refereed to as the \emph{Kummer sublocus} or the set of Kummer configurations of six lines. For a general member $Y$ of the Kummer locus the N\'eron-Severi lattice is given by Equation~(\ref{eqn:lattice}). Conversely, a K3 surface $\mathcal{X}$ with $\mathrm{NS}(\mathcal{X}) \cong L$ with $L$ in Equation~(\ref{eqn:lattice}) is a generic Jacobian Kummer surface and isomorphic to a double sextic surface branched on a Kummer configuration of six lines.
\section{Geometry of Period Domains for Lattice Polarized K3 surfaces}
\label{s-period_domains}
 In this section we review aspects of the period map for lattice polarized K3 surfaces. The corresponding periods are the solutions of  the Picard-Fuchs system. We then review Sasaki \& Yoshida's \cite{MR0996019} differential geometric interpretation of such equations, which equips the moduli space with a canonical integrable holomorphic conformal structure.
\subsection{Period Domains and the Period Map}
\label{ss-periods}
Let us start by reviewing how the transcendental lattice of a general member in a family of lattice polarized K3 surfaces determines a local deformation of complex structure. 
\par Let $L \overset{\jmath}{\hookrightarrow} \Lambda_{\mathrm{K3}}$ be a primitive sublattice of rank $\rho$ of the K3 lattice, and let $X$ be a general member in the family of $L$-polarized K3 surfaces, that is also equipped with a marking as in \S \ref{s-k3_kummer}. Let $N = L^\perp$ be the orthogonal complement in the K3 lattice.. Assume that $N$ is an even lattice of signature $(2,n)$ with $n \in \mathbb{N}$ so that
\begin{equation*}
	n + 2 = \mathrm{rank}(N) = \mathrm{rank}(\Lambda_{\mathrm{K3}}) - \mathrm{rank}(L) = 22 - \rho \, .
\end{equation*} 
\par Since $H^{2,0}(X) \perp L_{\mathbb{C}}$, the line $H^{2,0}(X) \subset N_{\mathbb{C}}$ is contained in the complexified transcendental lattice. If $\pi : N_{\mathbb{C}} - \{0\} \to \mathbb{P}(N_{\mathbb{C}})$ is the canonical projection, it follows from the Torelli theorem of Burns \& Rapoport for algebraic K3 surfaces \cite{ASENS_1975_4_8_2_235_0} that the marked K3 surface $X$ determines a unique point in $\mathbb{P}(N_\mathbb{C})$, called the period point, given by
\begin{equation}
	\label{eq-period_point}
	p_X = \pi\big(H^{2,0}(X)\big) \  \in \ \mathbb{P}(N_\mathbb{C}) = \mathbb{P}^{n+1} \, .
\end{equation}
Geometrically, the period point $p_X \in \mathbb{P}^{n+1}$ can be described as follows. Choosing a global nonvanishing section of the canonical bundle $\omega_X$ is equivalent to a nonvanishing holomorphic 2-form $\Omega$ with $[\Omega] \in H^{2,0}(X)$. Using Poincar\'e duality, we have $H_2(X,\mathbb{Z}) \cong N \oplus L$. We fix an ordered basis $\{\gamma_0,\dots,\gamma_{n+1},\gamma_{n+2}, \dots,\gamma_{21}\} \subset H_2(X,\mathbb{Z})$ of 2-cycles so that $N \cong \mathrm{span}_\mathbb{Z}\{\gamma_0,\dots,\gamma_{n+1}\}$ and $L \cong  \mathrm{span}_\mathbb{Z}\{\gamma_{n+2},\dots,\gamma_{21}\}$. It follows that the integral of $\Omega$ over any algebraic cycle $\gamma \in L$ vanishes, i.e., 
\begin{equation*}
	\oint_{\gamma} \Omega =0 \, ,
\end{equation*}
and the only nontrivial period integrals arise by integrating over transcendental cycles. In this way, the period point can also be obtained as
\begin{equation}
	\label{eq-period_integrals}
	p_X := \left[\,\oint_{\gamma_0} \! \! \Omega  : \,  \cdots \, :  \oint_{\gamma_{n+1}} \! \! \Omega  \right] \in \mathbb{P}^{n+1} \, .
\end{equation}
\par Since $[\omega \wedge \omega] \in H^{4,0}(X) = \{ 0\}$ we have $(\Omega,\Omega)= 0$, and it follows that $p_X$ lies on an $n$-dimensional quadric hypersurface in $\mathbb{P}^{n+1}$ which we will denote by $Q^n \subset \mathbb{P}^{n+1}$. The hypersurface $Q^n$ is defined as the projectivized locus of isotropic vectors in $N_{\mathbb{C}}$ with respect to the quadratic form $( \, \cdot \, , \, \cdot \, )$. Moreover, since 
\begin{equation}
	\label{eq-period_relation}
	\left(\Omega,\overline{\Omega}\right) > 0 \, ,
\end{equation}
it follows that $p_X \in \mathcal{D}_L$, where $\mathcal{D}_L$ is the relatively open subset of $Q^n$ defined by Equation (\ref{eq-period_relation}). Hence, $\mathcal{D}_L$ is a quasi-projective subvariety of $Q^n$, called the period domain. The period domain inherits more geometric structure than that of a quasi-projective variety. In fact, identifying $H^{2,0}(X) \subset N_{\mathbb{C}}$ with a positive definite oriented 2-plane in $N_{\mathbb{R}}$, it follows that $\mathcal{D}_L$ can be identified with the symmetric homogeneous space $\mathcal{D}_L \cong \mathrm{O}(2,n) / (\mathrm{SO}(2) \times \mathrm{O}(n))$ of positive oriented 2-planes in $N_{\mathbb{R}}$. We refer the interested reader to Aspinwall \cite[\S 2.5]{MR1479699} for an overview and a discussion of the connection with the moduli space of K\"ahler-Einstein metrics on a smooth K3 surface. Then $\mathcal{D}_L$ has two connected components which are exchanged under the induced involution by complex conjugation on forms. Each component is isomorphic to a bounded Hermitian symmetric domain of type $IV_{n}$. The coarse moduli space $\mathfrak{M}_L$ introduced in \S \ref{s-k3_kummer} then takes the form of $\mathfrak{M}_L = \mathcal{D}_L / \Gamma_L$, where $\Gamma_L$ is the kernel of the natural homomorphism $\mathrm{O}(L,\mathbb{Z}) \to A(L)$; see Dolgachev \cite[\S 3]{MR1420220} for a detailed discussion.
\par Next, we consider the family of marked $L$-polarized K3 surfaces $f_L \, \colon \, \mathfrak{X} \to M$ over a complex quasi-projective variety $M$ where $f_L$ is a proper smooth holomorphic map and each fiber over $t \in M$ is an $L$-polarized K3 surface $X_t$, that is also equipped with a basis $\{\gamma_i(t)\}_{i=0}^{21} \subset H_2(X_t ,\mathbb{Z})$ of locally constant 2-cycles on which to integrate a non-degenerate holomorphic 2-form $\Omega_t \in H^{2,0}(X_t)$; see Dolgachev \cite[Remark 3.4]{MR1420220} for more details. The holomorphic  map  
\begin{equation}
	\label{eq-period_map}
	\begin{aligned}
	p \; & \colon \quad M \to  \ \mathcal{D}_L \subset Q^n \, \\
	 t  \ \mapsto \ p_{X_t} & = \left[\,\oint_{\gamma_0(t)} \! \! \Omega_t  : \,  \cdots \, :  \oint_{\gamma_{n+1}(t)} \! \! \Omega_t  \right] \ \in \ \mathcal{D}_L \, ,
	\end{aligned} 
\end{equation}
is called the period map of the family $f_L \, \colon \, \mathfrak{X} \to M$. The family is called maximal if the period mapping is surjective. In this case we must that the family is not isotrivial and that $\dim(M) = n$, the latter by the local Torelli theorem of Burns \& Rapoport \cite{ASENS_1975_4_8_2_235_0} which guarantees that the period mapping is locally an isomorphism. In this way, the period map of a family encodes the local variation of the complex structure in the K3 fibers of the family $f_L$. 
\par By the global Torelli theorem (loc.~cit.), the period map also characterizes the variation of Hodge structure, that is, how Hodge groups $H^{2,0}(X_t),H^{1,1}(X_t),H^{0,2}(X_t)$ vary within $H^2(X_{t_0},\mathbb{C})$ (the latter can be described in purely topological terms) as $t \in M$ varies. Moreover, for $M = \mathfrak{M}_L$ the solutions of the Picard-Fuchs system determine a (multivalued) inverse of the canonical projection $\mathcal{D}_L \to \mathfrak{M}_L$ and provide a uniformization of the moduli space. 
\subsection{Linear PDEs modeled after quadric hypersurfaces}
\label{ss-sasaki_yoshida}
Inspired at least partially by the discussion above, Sasaki \& Yoshida \cite{MR0996019} sought to develop a \emph{geometric} theory of linear PDEs in $n$ variables of rank $n+2$, whose solutions provide a map of the equation manifold $M$ into a quadric hypersurface $Q^n \subset \mathbb{P}^{n+1}$. Such a theory had already been developed classically by Wilczynski \cite{MR0131232} for $n=1,2$, and Sasaki \& Yoshida successfully generalized Wilczynski's results for $n \geq 3$ by casting the equations in a canonical differential geometric form. We focus here on the case $n \geq 3$, which has a markedly different flavor than $n=1,2$. Sasaki \& Yoshida brought Wilczynski's results into the modern era for $n=2$ in \cite{MR0960834}. See Doran \cite{MR1779161,MR1738862} for aspects related to the $n=1$ case, and Clingher, Doran, and the first author \cite{MR3767270} for some aspects related to the $n=2$ case. The second author \cite[\S 6.2.3]{schultz_geometry_2021} connected some of these results with a particular $n=3$ case relevant to \S \ref{ss:double_sextics} that is briefly discussed at the end of \S \ref{ss-Kummer_PF}.
\par Let $M$ be a complex $n$-dimensional manifold with $n \geq 3$, equipped with a locally finite open covering by local holomorphic coordinate charts. Let $E$ be a second order holonomic system on $X$ in a single unknown of rank $n+2$, as defined by Sasaki \& Yoshida. In a local coordinate chart $(t^1,\dots,t^n)$, the system $E$ is then given by a second order system of linear homogeneous differential equations with polynomial coefficients (satisfying a suitable dimension condition in terms of D-modules theory) for a single unknown function $w=w(t^1,\dots,t^n)$ of rank $n+2$, whose solutions $w_0,\dots,w_{n+1}$ are unique up to scale. Therefore, it is natural to consider the following map determined by the solutions of the holonomic system:
\begin{equation}
	\label{eq-solution_map}
	\begin{aligned}
		p \; & \colon \quad M \to  \  \mathbb{P}^{n+1} \,, \\
		t  \ \mapsto \  p(t)  & = \left[\,w_0(t)\, : \,  \cdots \, : \, w_{n+1}(t) \, \right] \in \mathbb{P}^{n+1} \, .
	\end{aligned} 
\end{equation}
\par Sasaki \& Yoshida asked the following question: what is the geometry of $p(M) \subset \mathbb{P}^{n+1}$, and, conversely, what effect does imposing the existence of a geometric structure on $p(M)$ have on the solutions of $E$? Here, we will always assume that $p(M) \subseteq V \subset \mathbb{P}^{n+1}$, where $V$ is a smooth, closed complex hypersurface, $p(M)$ is relatively open in $V$, and $p$ is locally an isomorphism. By a celebrated result of Chow \cite{MR0033093}, $V$ is then a projective \emph{subvariety} of $\mathbb{P}^{n+1}$, and is given by the smooth locus of some algebraic equations. Hence, $p(M) \subseteq  V$ is a complex quasi-projective variety, so locally near a point $t \in M$, the solutions $w_0,\dots,w_{n+1}$ satisfy the same algebraic relations that also define $V \subset \mathbb{P}^{n+1}$ around $p(t) \in V$.
\par The simplest nontrivial algebraic relationship that can hold between solutions $w_0,\dots,w_{n+1}$ is a  quadratic relation. For example, the period integrals of lattice polarized K3 surfaces in \S \ref{ss-periods} always satisfy such quadratic relations. In general, for a second order holonomic system $E$ on $X$, we say that $E$ satisfies the \emph{quadric condition} if its solutions satisfy a quadratic relation. We can then take $V = Q^n \subset \mathbb{P}^{n+1}$ to be a quadratic hypersurface. Among Sasaki \& Yoshida's results are necessary and sufficient conditions \cite[Proposition 1.1]{MR0996019} to ensure that $p(M) \subseteq Q^n$. These conditions can be stated entirely in differential geometric terms on $M$. In fact, the relevant geometric structure to consider turns out to be the one of a \emph{holomorphic conformal structure} on the holomorphic tangent bundle $T_{\mathbb{C}}M$. The existence of such a structure for a complex manifold is quite restrictive \cite{MR0685426}. This is quite different from the situation for real (semi-)Riemannian manifolds that, as it is well known, always admit natural conformal structures.
\par We follow Sasaki \& Yoshida's presentation \cite[\S 1.1]{MR0996019}: let $h = (h_{ij}) \in S^2(\mathbb{C}^n)$ be the matrix representation of a fixed nondegenerate symmetric bilinear form on $\mathbb{C}^n$ with respect to some fixed basis, and $S^2(\mathbb{C}^n)$ denote the collection of all such forms. The conformal group $\mathrm{CO}(h)$ is defined to be
\begin{equation*}
	\mathrm{CO}(h)=\left\{\lambda a \mid a \in \mathrm{GL}(n, \mathbb{C}), a h a^T=h, \lambda \in \mathbb{C}^*\right\} \, , 
\end{equation*} 
which is the direct product of the multiplicative group $\mathbb{C}^\times$ and the complex orthogonal group $\mathrm{O}(h,\mathbb{C})$.
\par We say that $M$, a complex $n$-dimensional manifold with $n \ge 3$, admits a holomorphic conformal structure (HCS) if the structure group of the holomorphic tangent bundle reduces as a holomorphic bundle to the conformal group. This is equivalent to the condition that the complex frame bundle $F(M) \to M$ admits a holomorphic principal subbundle $P \to M$ with structure group $\mathrm{CO}(h)$ for some nondegenerate bilinear form $h$. As the bundle is topologically trivial, it is naturally associated with a holomorphic section $g\colon M \to F(M)/\mathrm{CO}(h)$, and hence, any HCS on $M$ corresponds to  a holomorphic, nondegenerate symmetric covariant tensor $g$ of rank-2, called the conformal metric. In a local coordinate chart, we write $g$ as
\begin{equation}
	\label{eq-conformal_metric}
	g = g_{ij} \, dt^i \odot dt^j \in S^2(M) \, .
\end{equation}
Here, $S^2(M)$ is the holomorphic bundle of symmetric rank-2 covariant tensors on ~$T_{\mathbb{C}}M$, $\odot$ is the symmetric tensor product, and we employ the summation convention. 
\par We call a HCS $g \in S^2(M)$ integrable if $g$ is locally conformally flat. It is then well that integrability of a metric can be detected from the vanishing of certain curvature tensors associated with $g$, namely the Weyl tensor $W(g)$ if $n \geq 4$, and the Cotton tensor $C(g)$ if $n=3$. Recall that these tensors are defined as follows: let $\Gamma_{i k}^j$ be the Christoffel symbols of $g$, i.e.,
\begin{equation}
	\label{eq-christoffel}
	\Gamma_{i k}^j=\frac{1}{2} g^{j l}\left(\partial_k g_{i l}+\partial_i g_{k l}-\partial_l g_{i k}\right) \,.
\end{equation}
The Christoffel symbols define a unique torsion free symmetric connection $\nabla$ on $T_\mathbb{C}M$ associated with $g$. The Riemannian curvature tensor $R^l{ }_{i j k}$ is defined as 
\begin{equation}
	\label{eq-curvature_tensor}
	R^l{ }_{i j k}=\partial_j \Gamma_{k i}^l-\partial_k \Gamma_{j i}^l+\Gamma_{j r}^l \Gamma_{k i}^r-\Gamma_{k r}^l \Gamma_{j i}^r \, .
\end{equation}
Moreover, the Ricci curvature $R_{ij}$ and the scalar curvature $R$ are given by the following contractions of Equation (\ref{eq-curvature_tensor}) with $g$:
\begin{equation}
	\label{eq-ricci_scalar}
	R_{i j}=R^l{ }_{i l j}, \quad R=g^{i j} R_{i j} \, .
\end{equation} 
Finally,  the Weyl tensor $W(g)$ is given by 
\begin{equation}
	\label{eq-weyl}
	W^j{ }_{i k l}=R^j{ }_{i k l}+S_{i k} \delta_l^j-S_{i l} \delta_k^j+g_{i k} g^{j m} S_{m l}-g_{i l} g^{j m} S_{m k} , 
\end{equation}
and the Cotton tensor $C(g)$ is
\begin{equation}
	\label{eq-cotton}
	C_{i j k}=\left(\nabla_i S(g)\right)_{j k}-\left(\nabla_j S(g)\right)_{i k} \, ,
\end{equation}
where $S(g)$ is the Schouten tensor
\begin{equation*}
	S_{i k}=\frac{1}{n-2}\left(R_{i k}-\frac{R}{2(n-1)} g_{i k}\right) \, .
\end{equation*}
\par Every quadric hypersurfaces $Q^n \subset \mathbb{P}^{n+1}$ naturally possesses an integrable HCS. By definition, $Q^n$ is the projectivized locus of isotropic vectors in $\mathbb{C}^{n+2}$ of a nondegenerate symmetric bilinear form $\tilde{h} \in S^2(\mathbb{C}^{n+2})$. Under a suitable automorphism of $\mathbb{C}^{n+2}$, $\tilde{h}$ can be brought into the form
\begin{equation*}
	\tilde{h}=\left(\begin{array}{rrr}
		0 & 0 & -1 \\
		0 & h & 0 \\
		-1 & 0 & 0
	\end{array}\right) \, ,
\end{equation*}
where $h \in S^2(\mathbb{C}^n)$ is nondegenerate. Identifying $\mathbb{P}^{n+1} = \mathbb{P}(X^0,\dots,X^{n+1})$,  $Q^n$ is then given by the homogeneous equation
\begin{equation}
	\label{eq-quadric_eqn}
	Q^n \; \colon \quad -2 X^0X^{n+1} + h_{ij}X^iX^j = 0 \, .
\end{equation}
One checks that we always have $p_0=[0 \, : \cdots \, : 1] \in Q^n$. The stabilizer group $\mathrm{O}(\tilde{h}) = \left\{g \in \mathrm{GL}(n+2,\mathbb{C}) \mid g \tilde{h} g^T = \tilde{h}\right\}$ acts transitively on $Q^n$, with an isotropy group $H$ at $p_0$, given by matrices of the form
\begin{equation}
	\label{eq-isotropy}
	\left(\begin{array}{lll}
		\lambda & 0 & 0 \\
		b & a & 0 \\
		\mu & c^T & \nu
	\end{array}\right) \quad \text{with} \ \lambda \nu=1, \, a h a^T=h, \, b=\lambda a h c, \, \mu=\lambda c^T h c / 2 ,
\end{equation}
where $\lambda, \mu, \nu \in \mathbb{C}^\times$, $b, c \in \mathbb{C}^n$, $a \in \mathrm{GL}(n,\mathbb{C})$. 
This is the standard presentation of the quadric hypersurface $Q^n = \mathrm{O}(\tilde{h}) / H$ as a homogeneous space and shows the reduction of structure group for the frame bundle $F(M)$ to $H$. 
Furthermore,  the linear isotropy representation at $p_0 \in Q^n$ has a kernel $N$, given by matrices of the form
\begin{equation}
	\label{eq-isotropy_kernel}
	\left(\begin{array}{ccc} 
		\pm 1 & 0 & 0 \\
		b & \pm \mathbb{I}_n & 0 \\
		\mu & c^T & \pm 1
	\end{array}\right)  .
\end{equation}
It follows that $H/N \cong \mathrm{CO}(h)$. Hence, there is a further reduction of the structure group of the frame bundle to $\mathrm{CO}(h)$. 
\par The corresponding HCS is called the canonical conformal structure of the quadric hypersurface $Q^n$, and it can be explicitly constructed as follows: we consider the constant tensor $\varphi = -2 dX^0 \odot dX^{n+1}+h_{ij}dX^i \odot dX^j$ on $\mathbb{C}^{n+2} - \{0\}$, and $\sigma$ a local section of the tautological bundle $\mathcal{O}_{\mathbb{P}^{n+1}}(-1)$ with projection $\pi \, \colon \, \mathbb{C}^{n+2} - \{0\} \to \mathbb{P}^{n+1}$.  We then obtain a conformal metric on $Q^n \subset \mathbb{P}^{n+1}$ by writing $\mathsf{g} = \sigma^*\varphi |_{Q^n}$. Although $\mathsf{g}$ depends on $\sigma$, changing the section results in rescaling by a nonvanishing holomorphic function, and hence the conformal class of $\mathsf{g}$ is well defined. It follows easily from the construction that $\mathsf{g}$ is also locally conformally flat. 
\par To summarize, for a holonomic system $E$ of rank $n+2$ over $M$ with $n=\dim{M} \geq ~3$, satisfying the quadric condition, the parameter space $M$ can be equipped with an integrable HCS $g \in S^2(M)$ by setting $g = p^*\mathsf{g}$, for $p$ given in Equation~(\ref{eq-solution_map}). Furthermore, the converse holds as well, as the following result shows:
\begin{theorem}[Theorem 2.5, \cite{MR0996019}]
	\label{thm-geometric_pdes}
	Let $M$ be an $n \geq 3$-dimensional complex manifold with an integrable HCS $g \in S^2(M)$. The system of linear differential equations, given by
	\begin{equation}
		\label{eq-quadric_system}
		g_{i j}\left(w_{k l}-\Gamma_{k l}^p w_p+\frac{1}{n-2} R_{k l} w\right)=g_{k l}\left(w_{i j}-\Gamma_{i j}^p w_p+\frac{1}{n-2} R_{i j} w\right) \, ,
	\end{equation}
is a holonomic system of rank $n+2$ satisfying the quadric condition. Here, subscripts on $w$ indicate partial derivatives, and $\Gamma^i_{jk}$, $R_{ij}$ are respectively the Christoffel symbols and Ricci curvature of $g$, given in Equations (\ref{eq-christoffel},\ref{eq-ricci_scalar}).
\end{theorem}
\par It follows then from the discussion in \S \ref{ss-periods} that the Picard-Fuchs system of any maximal family of $L$-polarized K3 surfaces with $\mathrm{rank}(N) \geq 5=3+2$ satisfies the quadric condition and can be represented in the form of Equation (\ref{eq-quadric_system}) for an integrable HCS on $M$ determined by the pairing $(\, \cdot \, , \, \cdot \,)$ on $N_{\mathbb{C}}$. In particular, the period domain $\mathcal{D}_L \subset Q^n$ naturally carries an integrable HCS via the inclusion map, and fundamental properties of the period map in Equation (\ref{eq-period_map}) from \S \ref{ss-periods}, such as the local Torelli theorem, impart $M$ with such an integrable HCS $g \in S^2(M)$ by pullback. Hence, the existence of such a $g \in S^2(M)$ completely characterizes the local variation of complex structure for such families. Let us see an example. 
\section{Holomorphic Conformal Geometry of the Kummer Sublocus}
\label{s-hcs}
In this section we compute the HCS of the Kummer sublocus of the parameter space of double sextic K3 surfaces from Equation (\ref{eqn:R}) introduced in \S  \ref{ss:double_sextics}. We then utilize Sasaki \& Yoshida's geometric theory from \S \ref{ss-sasaki_yoshida} to compute the Picard-Fuchs system of the family of generic Jacobian Kummer surfaces from \S \ref{ss:Kummers}. 
\subsection{Computation of the HCS}
\label{ss-HCS_computation}
We use a different, birational presentation for the family $\mathcal{X}$ in Equation~(\ref{eqn:MATeqn}), which is similar to the one studied by Hoyt \& Schwartz \cite{MR1877757}.  In fact, a  top-dimensional open subspace of $\mathcal{M}$  in Equation~(\ref{eqn:M6}) is given by elements $[ (a_{i j}) ] \in \mathcal{M}$ of the form
\begin{equation*}
\begin{pmatrix}
1 & 1  & 1  & 0 & 0 & -d \\
0 & 1  & 0 & -1 & -1 & -1 \\
0 & 0  & 1  & a & b & c 
\end{pmatrix}.
\end{equation*}
Setting $z_1=x, z_2=-t, z_3=-1$ in Equation~(\ref{eqn-six_lines}) we obtain a family of double sextic surfaces,  given by the equation
\begin{equation}
\label{eq-Hoyt_Schwartz_K3}
	X \;\colon \quad y^2=x(x-1)(x-t)(t-a)(t-b)(t-c-d x) \, ,
\end{equation}
over the four-dimensional parameter space 
\begin{equation}
\label{eq-Hoyt_Schwartz_param_space}
	M =\left\{(a, b, c, d) \in \mathbb{C}^4 \mid \ a \neq b, \  (c, d) \neq(a, 0),(b, 0),(0,1)\right\}  \,.
\end{equation}
This family of projective K3 surfaces was studied by the authors in \cite[\S 3]{MR4494119}.
As in \S \ref{ss-periods}, we can construct a family of marked $L$-polarized K3 surfaces $f_L \, \colon \, \mathfrak{X} \to M$ where the general fiber, obtained for generic parameters $(a, b, c, d) \in M$, is a smooth K3 surface $X$ in Equation~(\ref{eq-Hoyt_Schwartz_param_space}) with $\mathrm{NS}(X) \cong L$ where $L$ was given in Equation~(\ref{eqn:L}), and we have set $\Omega_{a,b,c,d} = dt \wedge dx/y$. 
\par It follows from \cite[Proposition 3.5]{MR4494119} that there is birational map $\psi \, \colon \, X \dashrightarrow \mathcal{X}$; in the affine coordinate system $z_1 =-1$, it is given by 
\begin{equation}
	\label{eq-K3_birational_map}
\psi \; \colon \quad 	z_2=\frac{x_3 t-1}{x_3 t-x_1}, \quad z_3=\frac{x\left(1-x_1\right)}{x_3 t-x_1}, \quad \tilde{z}_4=\frac{x_3\left(x_1-1\right)^2 y}{\left(x_3 t-x_1\right)^3},
\end{equation}
induces a birational map $\Psi \; \colon M \dashrightarrow \mathcal{M}$ between $M$ and $\mathcal{M}$, given by
\begin{equation}
\label{eq-induced_birational_map}
	\Psi \; \colon \quad x_1=\frac{a}{b}, \quad x_2=\frac{a-c}{b-c}, \quad x_3=\frac{1}{b}, \quad x_4=\frac{d}{b-c} .
\end{equation}
\par Matsumoto, Sasaki, \& Yoshida \cite[Equation (0.8.4)]{MR1136204} found an integrable holomorphic conformal structure $G = G_{ij}\, dx_i \odot dx_j \in S^2(\mathcal{M})$ on $\mathcal{M}$: they used the GKZ formalism \cite{MR1080980} for the period integrals of the family over $\mathcal{M}$ whose fibers are the projective K3 surfaces, determined by Equation~(\ref{eqn:MATeqn}). The  components $G_{ij}=G_{ji}$ of the HCS are as follows:
\begin{equation}
\label{eq-M_6_HCS}
	\begin{aligned}
	& G_{12}=\frac{x_4-x_3}{x_1-x_2}, \quad G_{13}=\frac{x_4-x_2}{x_1-x_3}, \quad G_{14}=1, \\
	& G_{23}=1, \quad G_{24}=\frac{x_3-x_1}{x_2-x_4}, \quad G_{34}=\frac{x_2-x_1}{x_3-x_4}, \\
	& G_{11}=\frac{x_2 x_3-x_4}{x_1\left(1-x_1\right)}-\frac{x_3\left(x_4-x_2\right)}{x_1\left(x_1-x_3\right)}-\frac{x_2\left(x_4-x_3\right)}{x_1\left(x_1-x_2\right)}, \\
	& G_{22}=\frac{x_1 x_4-x_3}{x_2\left(1-x_2\right)}-\frac{x_1\left(x_3-x_4\right)}{x_2\left(x_2-x_1\right)}-\frac{x_4\left(x_3-x_1\right)}{x_2\left(x_2-x_4\right)}, \\
	& G_{33}=\frac{x_1 x_4-x_2}{x_3\left(1-x_3\right)}-\frac{x_1\left(x_2-x_4\right)}{x_3\left(x_3-x_1\right)}-\frac{x_4\left(x_2-x_1\right)}{x_3\left(x_3-x_4\right)}, \\
	& G_{44}=\frac{x_2 x_3-x_1}{x_4\left(1-x_4\right)}-\frac{x_3\left(x_1-x_2\right)}{x_4\left(x_4-x_3\right)}-\frac{x_2\left(x_1-x_3\right)}{x_4\left(x_4-x_2\right)} \; .
	\end{aligned}
\end{equation}
Then $G$ is nondegenerate on $M$, and integrability is verified checking that the vanishing of the Weyl tensor, $W(G)=0$ from Equation (\ref{eq-weyl}). We then obtain an integrable holomorphic conformal structure $g$ on $M$ by pullback via $\Psi$, i.e., we set $g =~ \Psi^* G \in S^2(M)$.  The integrability of $g$ follows again by verifying the vanishing of the Weyl tensor, i.e., $W(g) = 0$.
\par The Kummer sublocus $M' $ within $M$ can be computed as the vanishing locus of the projective invariant $R$ in Equation~(\ref{eqn:R}). Notice that for $(a, b, c, d) \in M$ we have $a \neq b$; it follows from \cite[Proposition 3.12]{MR4494119} that the Kummer locus in $M$ is given by
\begin{equation}
 \label{eq-Kummer_sublocus}
 (ab - c) d - (b - c) (a - c) = 0 \,.
\end{equation}
It was shown in \cite[Proposition 3.14]{MR4494119}  that the lattice polarization of $f_L \, \colon \, \mathfrak{X} \to M$ upon restriction to $M'$ extends to the lattice
\begin{equation*}
	L'  \ \cong \  H \oplus D_8(-1) \oplus D_4(-1) \oplus A_3(-1)  \  \cong \ H \oplus D_7(-1) \oplus 2D_4(-1) \, .
\end{equation*}
We then construct a family of marked $L'$-polarized K3 surfaces $f_{L'} \, \colon \, \mathfrak{X}' \to M'$ where the fiber is the double sextic surface,  given by the equation
\begin{equation}
\label{eq-Hoyt_Schwartz_K3_2}
	X' \colon \quad y^2=x(x-1)(x-t)(t-a)(t-b)\big((ab-c)(t-c)-  (b - c) (a - c) x\big) \, ,
\end{equation}
and the family varies over the three-dimensional parameter space 
\begin{equation}
\label{eq-Kummer_param_space}
	M' =\left\{(a, b, c) \in \mathbb{C}^3 \mid \ a \neq b, \ a \neq c, \ b \neq c, \  c \neq ab \right\} .
\end{equation}
The general fiber, obtained for generic parameters $(a, b, c) \in M'$, is a smooth K3 surface $X'$ with $\mathrm{NS}(X') \cong L'$ of Picard number 17. A choice of holomorphic 2-form generating the period integrals on $X^{\prime}$ is $\Omega^{\prime}_{a,b,c} = dt \wedge dx/y$. As explained in \S \ref{s-k3_kummer}, $f_{L'}$ is a family of Jacobian Kummer surfaces.
\par Let $\imath \,\colon \, M' \hookrightarrow M$ be the inclusion map determined by Equation (\ref{eq-Kummer_sublocus}). A conformal metric on $M'$ is then given by
\begin{equation}
	\label{eq-Kummer_HCS}
		g' = \Big(-2bc^2(b-c)(b-1)(a-1)\Big)^{-1} \imath^*g \, \in S^2(M^{\prime}) \, .
\end{equation}
The components $g'_{ij} = g'_{ji}$, with $i, j \in \{ a, b, c\}$ are as follows:
\begin{equation}
	\label{eq-Kummer_HCS_components}
	\begin{aligned}
		&g'_{aa} =1, \quad g'_{ab} = 0 = g'_{cc}, \quad g'_{ac} = -\frac{a}{2c}, \\
		&g'_{bb} = - \frac{a(a-c)(a-1)}{b(b-c)(b-1)}, \quad g'_{bc} = \frac{a(a-c)(a-1)}{2c(b-c)(b-1)} \, . \\
	\end{aligned}
\end{equation}
Here we have chosen the additional holomorphic conformal scaling factor in front of $\imath^*g$ to simplify the metric. It follows that $g'$ defines a nondegenerate holomorphic conformal structure on $M'$. First, one checks that the determinant of the metric is nonzero everywhere; this can be seen from Equation (\ref{eq-Kummer_param_space}).  Second, one verifies the integrability of $g'$ by checking that the Cotton tensor $C(g')=0$ vanishes, using Equation (\ref{eq-cotton}).
\subsection{Picard-Fuchs Equations for Generic Kummer Surfaces}
\label{ss-Kummer_PF}
We can now use the results of Theorem \ref{thm-geometric_pdes} with $n=3$ and $g^{\prime} \in S^2(M^{\prime})$ from Equation (\ref{eq-Kummer_HCS_components}) to characterize explicitly the Picard-Fuchs system for the family of Jacobian Kummer surfaces $f_{L'} \, \colon \, \mathfrak{X}' \to M'$. To the best of our knowledge, such a system has not been determined before. 
\begin{theorem}
\label{thm-Kummer_PF}
The system of linear differential equations, given by
	\begin{equation}
		\label{eq-Kummer_PF_eqs}
		g'_{i j}\left(w_{kl}-\Gamma_{kl}^r w_r+ R_{kl} w\right)=g'_{kl}\left(w_{i j}-\Gamma_{i j}^r w_r+ R_{i j} w\right) \quad \text{for $i, j, k, l, r \in \{ a, b, c\}$},
	\end{equation} 	
is a rank-5 holonomic system in the variables $a, b, c$. It is the Picard-Fuchs system annihilating the period integrals for the family $f_{L'} \,\colon \, \mathfrak{X}' \to~ M'$, and any other representation of the Picard-Fuchs system is projectively equivalent to Equation (\ref{eq-Kummer_PF_eqs}).  Here, $g' \in S^2(M')$ is the integrable holomorphic conformal structure in Equation (\ref{eq-Kummer_HCS_components}), and $\Gamma^r_{ij}$, $R_{ij}$ are the Christoffel connection and Ricci curvature of $g'$, respectively.
\end{theorem}
\par The Picard-Fuchs system in Equation (\ref{eq-Kummer_PF_eqs}) necessarily satisfies the quadric condition by Theorem \ref{thm-geometric_pdes} since the HCS is integrable. Such quadratic relationships, when attached to Picard-Fuchs equations, always contain intricate geometric and arithmetic information about the K3 surface itself. Sasaki \& Yoshida \cite[Equation 3.1]{MR0996019} determined a closely related system of rank-5, acting as the uniformizing equations for the Siegel modular threefold $\mathcal{A}_2[2] = \mathbb{H}_2 / \Gamma(2)$. Here, $\Gamma(2) \subset \mathrm{Sp}(4,\mathbb{Z})$ is the Siegel modular group of level-two. It follows from \cite[\S 0.19]{MR1136204} that this system is the Picard-Fuchs system for Kummer surfaces associated with the Jacobian of smooth genus-2 curves $C$ with level-2 structure as in Equation~(\ref{eqn:Rosenhain}). Along with Hara \cite{MR1040172}, they showed that their system decomposes as the exterior product of two copies of Lauricella's $F_D$ system, each of which represents precisely the Picard-Fuchs system of such a genus-2 curve $C$.
\par In \cite[Theorem 6.2.89]{schultz_geometry_2021}, the second author showed that the rank-5 uniformizing system for $\mathcal{A}_2[2]$ is closely related to the Picard-Fuchs system of the subfamily $X'' =~ X_{a,b,c,0}$ of Equation (\ref{eq-Hoyt_Schwartz_K3}). This family is also known as the twisted Legendre pencil, and it was originally studied by Hoyt \cite{MR1013162}. The family does \emph{not} represent Kummer surfaces of a principally polarized abelian surface, since Equation (\ref{eq-Kummer_sublocus}) does not hold. Rather, Equation (\ref{eq-Hoyt_Schwartz_K3}) with $d=0$ represents a six line configuration in $\mathbb{P}^2$ where three lines intersect in a single point \cite{clingher_geometry_2017}, and $X''$ gives the relative Jacobian fibration of a Kummer surface of a $(1, 2)$-polarized abelian surface; see \cite{MR4484238}. Then under a change of variables determined by Braeger, Sung, and the first author \cite[Eq.~(3.11)]{MR4099481}, the second author showed (loc. cit.) that the Picard-Fuchs system for $X''$ also decomposes into the exterior product of Lauricella $F_D$ systems.
\subsection{Future directions}
At present it is not known how a multivariate system like the one in Equation (\ref{eq-quadric_system}) decomposes in general when imposing the quadric condition. It would be interesting to determine precisely how the Picard-Fuchs system in Equation (\ref{eq-Kummer_PF_eqs}) decomposes under the quadric condition.  One might  also ask how the  Picard-Fuchs system decomposes further when a Kummer configuration of six lines in the projective plane attains certain special curves; essential information of the principally polarized abelian surface $A$ are reflected in these configurations of six lines, used to construct $\mathrm{Kum}(A)$ as a double sextic surface. For example, a geometric characterizations of the Humbert surfaces, within the moduli space $\mathcal{A}_2$ of principally polarized abelian surfaces, in terms of the presence of certain curves on the associated Kummer plane was given in \cite{MR1953527}.  These questions  will be discussed in forthcoming works by the authors.
\newpage
\bibliographystyle{abbrv}
\bibliography{MalmendierSchultz_HCS}

\begin{thebibliography}{10}

\bibitem{MR1479699}
P.~S. Aspinwall.
\newblock {$K3$} surfaces and string duality.
\newblock In {\em Fields, strings and duality ({B}oulder, {CO}, 1996)}, pages
  421--540. World Sci. Publ., River Edge, NJ, 1997.

\bibitem{MR2062673}
C.~Birkenhake and H.~Lange.
\newblock {\em Complex abelian varieties}, volume 302 of {\em Grundlehren der
  mathematischen Wissenschaften [Fundamental Principles of Mathematical
  Sciences]}.
\newblock Springer-Verlag, Berlin, second edition, 2004.

\bibitem{MR1953527}
C.~Birkenhake and H.~Wilhelm.
\newblock Humbert surfaces and the {K}ummer plane.
\newblock {\em Trans. Amer. Math. Soc.}, 355(5):1819--1841, 2003.

\bibitem{MR4099481}
N.~Braeger, A.~Malmendier, and Y.~Sung.
\newblock Kummer sandwiches and {G}reene-{P}lesser construction.
\newblock {\em J. Geom. Phys.}, 154:103718, 18, 2020.

\bibitem{ASENS_1975_4_8_2_235_0}
D.~Burns and M.~Rapoport.
\newblock On the {Torelli} problem for {K}\"ahlerian ${K}3$ surfaces.
\newblock {\em Annales scientifiques de l'\'Ecole Normale Sup\'erieure}, Ser.
  4, 8(2):235--273, 1975.

\bibitem{MR0033093}
W.-L. Chow.
\newblock On compact complex analytic varieties.
\newblock {\em Amer. J. Math.}, 71:893--914, 1949.

\bibitem{MR3767270}
A.~Clingher, C.~F. Doran, and A.~Malmendier.
\newblock Special function identities from superelliptic {K}ummer varieties.
\newblock {\em Asian J. Math.}, 21(5):909--951, 2017.

\bibitem{clingher_geometry_2017}
A.~Clingher and A.~Malmendier.
\newblock On the geometry of (1,2)-polarized {Kummer} surfaces, Apr. 2017.
\newblock arXiv:1704.04884 [hep-th].

\bibitem{MR4015343}
A.~Clingher, A.~Malmendier, and T.~Shaska.
\newblock Six line configurations and string dualities.
\newblock {\em Comm. Math. Phys.}, 371(1):159--196, 2019.

\bibitem{MR4484238}
A.~Clingher, A.~Malmendier, and T.~Shaska.
\newblock On isogenies among certain {A}belian surfaces.
\newblock {\em Michigan Math. J.}, 71(2):227--269, 2022.

\bibitem{MR1007155}
I.~Dolgachev and D.~Ortland.
\newblock Point sets in projective spaces and theta functions.
\newblock {\em Ast\'{e}risque}, (165):210 pp. (1989), 1988.

\bibitem{MR1420220}
I.~V. Dolgachev.
\newblock Mirror symmetry for lattice polarized {$K3$} surfaces.
\newblock volume~81, pages 2599--2630. 1996.
\newblock Algebraic geometry, 4.

\bibitem{MR1779161}
C.~F. Doran.
\newblock Picard-{F}uchs uniformization and modularity of the mirror map.
\newblock {\em Comm. Math. Phys.}, 212(3):625--647, 2000.

\bibitem{MR1738862}
C.~F. Doran.
\newblock Picard-{F}uchs uniformization: modularity of the mirror map and
  mirror-moonshine.
\newblock In {\em The arithmetic and geometry of algebraic cycles ({B}anff,
  {AB}, 1998)}, volume~24 of {\em CRM Proc. Lecture Notes}, pages 257--281.
  Amer. Math. Soc., Providence, RI, 2000.

\bibitem{MR1080980}
I.~M. Gel'fand, M.~M. Kapranov, and A.~V. Zelevinsky.
\newblock Generalized {E}uler integrals and {$A$}-hypergeometric functions.
\newblock {\em Adv. Math.}, 84(2):255--271, 1990.

\bibitem{MR229641}
P.~A. Griffiths.
\newblock Periods of integrals on algebraic manifolds. {I}. {C}onstruction and
  properties of the modular varieties.
\newblock {\em Amer. J. Math.}, 90:568--626, 1968.

\bibitem{MR258824}
P.~A. Griffiths.
\newblock Periods of integrals on algebraic manifolds: {S}ummary of main
  results and discussion of open problems.
\newblock {\em Bull. Amer. Math. Soc.}, 76:228--296, 1970.

\bibitem{MR1040172}
M.~Hara, T.~Sasaki, and M.~Yoshida.
\newblock Tensor products of linear differential equations---a study of
  exterior products of hypergeometric equations.
\newblock {\em Funkcial. Ekvac.}, 32(3):453--477, 1989.

\bibitem{MR1013162}
W.~L. Hoyt.
\newblock On twisted {L}egendre equations and {K}ummer surfaces.
\newblock In {\em Theta functions---{B}owdoin 1987, {P}art 1 ({B}runswick,
  {ME}, 1987)}, volume 49, Part 1 of {\em Proc. Sympos. Pure Math.}, pages
  695--707. Amer. Math. Soc., Providence, RI, 1989.

\bibitem{MR1877757}
W.~L. Hoyt and C.~F. Schwartz.
\newblock Yoshida surfaces with {P}icard number {$\rho\geq 17$}.
\newblock In {\em Proceedings on {M}oonshine and related topics
  ({M}ontr\'{e}al, {QC}, 1999)}, volume~30 of {\em CRM Proc. Lecture Notes},
  pages 71--78. Amer. Math. Soc., Providence, RI, 2001.

\bibitem{MR3586372}
D.~Huybrechts.
\newblock {\em Lectures on {K}3 surfaces}, volume 158 of {\em Cambridge Studies
  in Advanced Mathematics}.
\newblock Cambridge University Press, Cambridge, 2016.

\bibitem{MR0355886}
S.~Kobayashi.
\newblock {\em Transformation groups in differential geometry}, volume Band 70
  of {\em Ergebnisse der Mathematik und ihrer Grenzgebiete [Results in
  Mathematics and Related Areas]}.
\newblock Springer-Verlag, New York-Heidelberg, 1972.

\bibitem{MR0685426}
S.~Kobayashi and T.~Ochiai.
\newblock Holomorphic structures modeled after hyperquadrics.
\newblock {\em Tohoku Math. J. (2)}, 34(4):587--629, 1982.

\bibitem{MR0031310}
N.~H. Kuiper.
\newblock On conformally-flat spaces in the large.
\newblock {\em Ann. of Math. (2)}, 50:916--924, 1949.

\bibitem{MR3263663}
A.~Kumar.
\newblock Elliptic fibrations on a generic {J}acobian {K}ummer surface.
\newblock {\em J. Algebraic Geom.}, 23(4):599--667, 2014.

\bibitem{MR2409557}
M.~Kuwata and T.~Shioda.
\newblock Elliptic parameters and defining equations for elliptic fibrations on
  a {K}ummer surface.
\newblock In {\em Algebraic geometry in {E}ast {A}sia---{H}anoi 2005},
  volume~50 of {\em Adv. Stud. Pure Math.}, pages 177--215. Math. Soc. Japan,
  Tokyo, 2008.

\bibitem{MR4494119}
A.~Malmendier and M.~T. Schultz.
\newblock On the mixed-twist construction and monodromy of associated
  {P}icard-{F}uchs systems.
\newblock {\em Commun. Number Theory Phys.}, 16(3):459--513, 2022.

\bibitem{MR1103969}
K.~Matsumoto.
\newblock Theta functions on the classical bounded symmetric domain of type
  {${\rm I}_{2,2}$}.
\newblock {\em Proc. Japan Acad. Ser. A Math. Sci.}, 67(1):1--5, 1991.

\bibitem{MR973860}
K.~Matsumoto, T.~Sasaki, and M.~Yoshida.
\newblock The period map of a {$4$}-parameter family of {$K3$} surfaces and the
  {A}omoto-{G}elfand hypergeometric function of type {$(3,6)$}.
\newblock {\em Proc. Japan Acad. Ser. A Math. Sci.}, 64(8):307--310, 1988.

\bibitem{MR1073363}
K.~Matsumoto, T.~Sasaki, and M.~Yoshida.
\newblock The {A}omoto-{G}elfand hypergeometric function and period mappings.
\newblock {\em S\={u}gaku}, 41(3):258--263, 1989.

\bibitem{MR1136204}
K.~Matsumoto, T.~Sasaki, and M.~Yoshida.
\newblock The monodromy of the period map of a {$4$}-parameter family of {$K3$}
  surfaces and the hypergeometric function of type {$(3,6)$}.
\newblock {\em Internat. J. Math.}, 3(1):164, 1992.

\bibitem{MR2514037}
D.~Mumford.
\newblock {\em Abelian varieties}, volume~5 of {\em Tata Institute of
  Fundamental Research Studies in Mathematics}.
\newblock Tata Institute of Fundamental Research, Bombay; by Hindustan Book
  Agency, New Delhi, 2008.
\newblock With appendices by C. P. Ramanujam and Yuri Manin, Corrected reprint
  of the second (1974) edition.

\bibitem{MR0429917}
V.~V. Nikulin.
\newblock Kummer surfaces.
\newblock {\em Izv. Akad. Nauk SSSR Ser. Mat.}, 39(2):278--293, 471, 1975.

\bibitem{MR0960834}
T.~Sasaki and M.~Yoshida.
\newblock Linear differential equations in two variables of rank four. {I}.
\newblock {\em Math. Ann.}, 282(1):69--93, 1988.

\bibitem{MR0996019}
T.~Sasaki and M.~Yoshida.
\newblock Linear differential equations modeled after hyperquadrics.
\newblock {\em Tohoku Math. J. (2)}, 41(2):321--348, 1989.

\bibitem{schultz_geometry_2021}
M.~T. Schultz.
\newblock {\em On the {Geometry} of the {Moduli} {Space} of {Certain} {Lattice}
  {Polarized} {K3} {Surfaces} and {Their} {Picard}-{Fuchs} {Operators}}.
\newblock PhD thesis, Utah State University, 2021.

\bibitem{MR0131232}
E.~J. Wilczynski.
\newblock {\em Projective differential geometry of curves and ruled surfaces}.
\newblock Chelsea Publishing Co., New York, 1962.

\end{thebibliography}
\end{document}